\documentclass[conference]{IEEEtran}
\usepackage{amsmath}
\usepackage{amssymb}
\usepackage{graphicx}

\ifCLASSINFOpdf
\else
\fi
\begin{document}
\title{A Comparison of Automatic Differentiation and Continuous Sensitivity
Analysis for Derivatives of Differential Equation Solutions}

\makeatletter
\newcommand{\linebreakand}{%
  \end{@IEEEauthorhalign}
  \hfill\mbox{}\par
  \mbox{}\hfill\begin{@IEEEauthorhalign}
}
\makeatother

\author{\IEEEauthorblockN{Yingbo Ma}
\IEEEauthorblockA{Julia Computing\\
Cambridge, Massachusetts, USA\\
Email: yingbo.ma@juliacomputing.com}
\and
\IEEEauthorblockN{Vaibhav Dixit}
\IEEEauthorblockA{Julia Computing\\
Pumas-AI\\
Bangalore, India\\
Email: vaibhav.dixit@juliacomputing.com}
\and
\IEEEauthorblockN{Michael J Innes}
\IEEEauthorblockA{Edinburgh, UK\\
Email: mike.j.innes@gmail.com}
\linebreakand 
\IEEEauthorblockN{Xingjian Guo}
\IEEEauthorblockA{New York University\\
New York City, New York, USA\\
Email: xg703@nyu.edu}
\and
\IEEEauthorblockN{Chris Rackauckas}
\IEEEauthorblockA{Massachusetts Institute of Technology\\
Julia Computing\\
Pumas-AI\\
Cambridge, Massachusetts, USA\\
Email: chris.rackauckas@juliacomputing.com}}

\maketitle

\begin{abstract}
    Derivatives of differential equation solutions are commonly for parameter estimation, fitting neural differential equations, and as model diagnostics. However, with a litany of choices and a Cartesian product of potential methods, it can be difficult for practitioners to understand which method is likely to be the most effective on their particular application. In this manuscript we investigate the performance characteristics of Discrete Local Sensitivity Analysis implemented via Automatic Differentiation (DSAAD) against continuous adjoint sensitivity analysis. Non-stiff and stiff biological and pharmacometric models, including a PDE discretization, are used to quantify the performance of sensitivity analysis methods. Our benchmarks show that on small stiff and non-stiff systems of ODEs (approximately $<100$ parameters+ODEs), forward-mode DSAAD is more efficient than both reverse-mode and continuous forward/adjoint sensitivity analysis. The scalability of continuous adjoint methods is shown to be more efficient than discrete adjoints and forward methods after crossing this size range. These comparative studies demonstrate a trade-off between memory usage and performance in the continuous adjoint methods that should be considered when choosing the technique, while numerically unstable backsolve techniques from the machine learning literature are demonstrated as unsuitable for most scientific models. The performance of adjoint methods is shown to be heavily tied to the reverse-mode AD method used for the vector-Jacobian product calculations, with tape-based AD methods shown to be 2 orders of magnitude slower on nonlinear partial differential equations than static AD techniques. In addition, these results demonstrate the out-of-the-box applicability of DSAAD to differential-algebraic equations, delay differential equations, and hybrid differential equation systems where the event timing and effects are dependent on model parameters, showcasing an ease of implementation advantage for DSAAD approaches. Together, these benchmarks provide a guide to help practitioners to quickly identify the best mixture of continuous sensitivities and automatic differentiation for their applications.
\end{abstract}


% \fi
\IEEEpeerreviewmaketitle

\section{Introduction\label{sec:Introduction}}

In the literature of differential equations, local sensitivity analysis is the practice of calculating derivatives
to a differential equation's solution with respect to model parameters.
For an ordinary differential equation (ODE) of the form
\begin{equation}
\dot{u}=f(u,p,t),\label{eq:ode}
\end{equation}
where $f$ is the derivative function and $p$ are the model parameters, the sensitivity of the state
vector $u$ with respect to model parameter $p_{i}$ at time $t$
is defined as $\frac{\partial u(p,t)}{\partial p_{i}}$. These sensitivities have many applications. For example, they can be directly utilized in fields such as biological modeling to identify
parameters of interest for tuning and experimentation \cite{sommer_numerical_nodate}. 
Recent studies have utilized these sensitivities as part for training neural networks
associated with the ODEs \cite{chen_neural_2018,grathwohl_ffjord:_2018}.
In addition, these sensitivities are
indirectly utilized in many disciplines for parameter estimation of dynamical models.
Parameter estimation is the problem of finding parameters $p$ such
that a cost function $C(p)$ is minimized (usually some fit against
data) \cite{peifer_parameter_2007,hamilton_parameter_2011,zhen_parameter_nodate,zimmer_parameter_2013,sommer_numerical_nodate,friberg_model_2002,steiert_experimental_2012}.
Gradient-based optimization methods require
the computation of gradients of $C(p)$. By the chain rule, $\frac{dC}{dp}$ requires the calculation of $\frac{du(t_i)}{dp}$ which are the model sensitivities. Given the high computational
cost of parameter estimation due to the number of repeated numerical
solutions which are required, efficient and accurate computation of
model sensitivities is an important part of differential equation
solver software.

The simplest way to calculate model sensitivities is to utilize numerical
differentiation which is given by the formula
\begin{equation}
\frac{\partial u(t)}{\partial p_{i}}=\frac{u(p+\Delta p_{i},t)-u(p,t)}{\Delta p_{i}}+\mathcal{O}(\Delta p_{i}),\label{eq:numdiff}
\end{equation}
where $p+\Delta p_{i}$ means adding $\Delta p_{i}$ to only the $i$th
component of $p$. However, this method is not efficient (it requires
two numerical ODE solutions for each parameter $i$) and it is prone
to numerical error. If $\Delta p_{i}$ is chosen too large, then the
error term of the approximation is large. In contrast, if $\Delta p_{i}$
is chosen too small, then calculations may exhibit floating point
cancellation which increases the error \cite{burden_numerical_2011}.

To alleviate these issues, many differential equation solver softwares
implement a form of sensitivity calculation called continuous local
sensitivity analysis (CSA) \cite{zhang_fatode:_2014,alan_c._hindmarsh_sundials:_2005}.
Forward-mode continuous sensitivity analysis calculates the model sensitivities
by extending the ODE system to include the equations:
\begin{equation}
\frac{d}{dt}\left(\frac{\partial u}{\partial p_{i}}\right)=\frac{\partial f}{\partial u}\frac{\partial u}{\partial p_{i}}+\frac{\partial f}{\partial p_{i}}\label{eq:clsa}
\end{equation}
where $\frac{\partial f}{\partial u}$ is the Jacobian of the derivative
function $f$ with respect to the current state, and $\frac{\partial f}{\partial p_{i}}$ 
is the derivative of the derivative function with respect to the $i$th
parameter. We note that $\frac{\partial f}{\partial u} v$ is equivalent to the
directional derivative in the direction of $v$, which can thus be calculated Jacobian-free
via: 
\begin{equation}
\frac{\partial f}{\partial u} v \approx \frac{f(u+\epsilon v,p,t)-f(u,p,t)}{\epsilon}
\end{equation}
or by equivalently pre-seeding forward-mode automatic differentiation with $v$ in the dual space.
Since the sensitivity equations for each $i$ are dependent on the
current state $u$, these ODEs must be solved simultaneously with
the ODE system $u'=f$. By solving this expanded system, one can ensure
that the sensitivities are computed to the same error tolerance as
the original ODE terms, and only a single numerical ODE solver call is required.

However, since the number of ODEs in this system now scales proportionally
with the number of parameters, the scaling of the aforementioned method is
$\mathcal{O}(np)$ for $n$ ODEs and $p$ parameters, forward-mode CSA is not
practical for a large number of parameters. Instead, for these cases
continuous adjoint sensitivity analysis (CASA) is utilized. This methodology
is defined to directly compute the gradient of a cost function of the solution with
scaling $\mathcal{O}(n+p)$.
Given a cost function $c$ on the ODE solution which is evaluated as discrete time points (such as an $L^2$ loss)

\begin{equation}
C(u(p), p) = \sum_i c(u(p,t_i), p)
\end{equation}
this is done by solving a backwards ODE known as adjoint problem

\begin{equation}
\frac{d\lambda'}{dt} = - \lambda' \frac{\partial f(u(t),p,t)}{\partial u}
\end{equation}
where at every time point $t_i$, this backwards ODE is perturbed by $\frac{\partial c(u(p,t_i), p)}{\partial u}$. Note $u(t)$ is generated by a forward solution. The gradient of the cost function is then given by the integral:

\begin{align}
\frac{dC}{dp} = &\lambda'(t_0) \frac{\partial f(u(t_0),p,t_0)}{\partial u} + \\ \nonumber
&\sum_i \int_{t_i}^{t_{i+1}} \lambda' \frac{\partial f(u(t),p,t)}{\partial p} + \frac{\partial c(u(p,t_i), p)}{\partial p} dt
\end{align}
This integral may be solved via quadrature on a continuous
solution of the adjoint equation or by appending the
quadrature variables to perform the integration as part
of the backsolve. The former can better utilize quadrature
points to reduce the number of evaluation points required
to reach a tolerance, but requires enough memory to store
a continuous extension of the backpass. We note that the backsolve
technique is numerically unstable, which will be noted in benchmarks
where it produces a divergent calculation of the gradient. We note that, 
similarly to forward-mode,
the crucial term $\lambda'(t) \frac{\partial f(u(t),p,t)}{\partial p}$ can be
computed without building the full Jacobian. However, in this case this computation
$v'J$ cannot be matrix-free computed via numerical or forward-mode differentiation,
but instead is the primitive of reverse-mode automatic differentiation. Thus efficient
methods for CASA necessarily mix a reverse-mode AD into the generated adjoint pass,
and thus we will test the Cartesian product of the various choices.

In contrast to CSA methods, discrete sensitivity
analysis calculates model sensitivities by directly differentiating the
numerical method's steps \cite{zhang_fatode:_2014}. However, this approach
requires specialized implementations of the first order ODE solvers to
propagate said derivatives. Instead, one can achieve the same end by using automatic differentiation
(AD) on a solver implemented entirely in a language with pervasive
AD, also known as a differentiable programming approach \cite{rackauckas2020generalized}. Section \ref{sec:Discrete-Sensitivity-Analysis} introduces the discrete sensitivity analysis
through AD (DSAAD) approach via type-specialization on a generic algorithm.
Section \ref{sec:Discrete-and-Continuous} compares the performance of DSAAD
against continuous sensitivity analysis and numerical differentiation
approaches and shows that DSAAD consistently performs well on the tested models. 
Section \ref{sec:Generalizes} describes limitations of the
continuous sensitivity analysis approach and describes how these cases are
automatically handled in the case of DSAAD. Together, this manuscript shows that
the ability to utilize AD directly on numerical integrator can be
advantageous to existing approaches for the calculation of model sensitivities, while at times can be disadvantageous in terms of scaling performance.

Currently, continuous sensitivity techniques are commonly used throughout
software such as SUNDIALS \cite{alan_c._hindmarsh_sundials:_2005}, while some software like FATODE \cite{zhang_fatode:_2014} allow for
discrete sensitivity analysis. No previous manuscript performs
comprehensive benchmarks on the Cartesian product of 
discrete/continuous forward/adjoint sensitivity analysis with the various automatic differentiation
modes (tape vs static). This study uses DiffEqSensitivity.jl \cite{rackauckas2020universal}, the sensitivity analysis extension to DifferentialEquations.jl and the first comprehensive package which includes all mentioned differentiation choices, to do such a full comparison of the space. The results showcase the efficiency gains provided by 
discrete sensitivity analysis on sufficiently small models, and
establishes a heuristic range ($\approx$ 30-100 ODEs) at which the scalability
of continuous adjoint sensitivity analysis overcomes the low overhead
of DSAAD. This study can thus can be a central to helping all users of ODE solvers to choose the method that will be effective on their specific problem.

\section{Discrete Sensitivity Analysis via Automatic Differentiation (DSAAD) \label{sec:Discrete-Sensitivity-Analysis}}

The core feature of the Julia programming language is multiple dispatch \cite{bezanson_julia:_2017}. It allows a function to compile to different outputs dependent
on the types of the inputs, effectively allowing choices of input
types to trigger forms of code generation. ForwardDiff.jl provides
a Dual number type which performs automatic differentiation on differentiable
programs by simultaneously propagating a derivative along with the
computed value on atomic (addition, multiplication, etc.) and standard
mathematical ($\sin$, $\exp$, etc.) function calls \cite{revels_forward-mode_2016}.
By using
the chain rule during the propagation, any function which is composed
of differentiable calls is also differentiable by the methodology.
Such a program is known as a differentiable program. Since this exactly
differentiates the atomics, the numerical error associated with this
method is similar to the standard evaluation of the function, effectively
alleviating the errors seen in numerical differentiation. In addition,
this method calculates derivatives simultaneously with the function's
evaluation, making it a good candidate for fast discrete sensitivity
analysis. It can be thought of as the analogue to continuous forward-mode
sensitivity analysis on the set of differentiable programs. Similarly,
reverse-mode AD from packages like ReverseDiff.jl \cite{revels_reversediff.jl_nodate} and Flux.jl \cite{innes_flux:_2018}
use Tracker numerical types which builds a tape
of the operations and utilizes the reverse of the chain rule to ``backpropogate''
derivatives and directly calculate the gradient of some cost function.
This implementation of AD can be thought of as the analogue of CASA on the
set of differentiable programs.

The DifferentialEquations.jl package provides many integration routines
which were developed in native Julia \cite{christopher_rackauckas_differentialequations.jl_2017}. These methods
are type-generic,
meaning they utilize the numeric and array types that are supplied
by the user. Thus these ODE solvers serve as a generic template whose
internal operations can be modified by external packages via dispatch.
When a generic DifferentialEquations.jl
integrator is called with a Dual number type for the initial condition,
the combination of these two programs results in a program which performs
discrete sensitivity analysis. This combination is what we define
as forward-mode DSAAD, and this combination with reverse-mode ReverseDiff.jl
AD Tracker types is reverse-mode DSAAD.

To test the correctness of the DSAAD and continuous sensitivity methods, we check the outputted
sensitivities on four models. For our tests we utilize nonlinear reaction
models which are representative of those found in biological and pharmacological
applications where these techniques are commonly used \cite{hiroaki_kitano_computational_nodate,danhof_mechanism-based_2008}. The models are:

\begin{enumerate}
    \item The non-stiff Lotka-Volerra equations (LV).
    \item An $N \times N$ finite difference discretization of the two-dimensional stiff Brusselator reaction-diffusion PDE (BRUSS).
    \item A stiff pollution model (POLLU).
    \item A non-stiff pharmacokinetic/pharmacodynamic system (PK/PD).
\end{enumerate}

These cover stiff and non-stiff ODEs, large systems and small systems, 
and include a PDE discretization with a dimension $N$ for testing 
the scaling of the methodologies. Details of the four models are 
presented in the Appendix.

Figure \ref{fig:ode_sens} shows the output of the first two models' 
sensitivities that
are computed by the DSAAD method compared to CSA. The two methods align in
their model sensitivity calculations, demonstrating that the application of AD
on the generic ODE solver does produce correct output sensitivities. From these
tests we note that the differences in the model sensitivities between the two
methods had a maximum norm of $1.14\times10^{-5}$ and $3.1\times10^{-4}$ which is
roughly the chosen tolerance of the numerical integration. These results were
again confirmed at a difference of $1\times10^{-12}$ when ran with sufficiently
low ODE solver tolerance.

\begin{figure}
    \begin{center}
        \includegraphics[width=\linewidth]{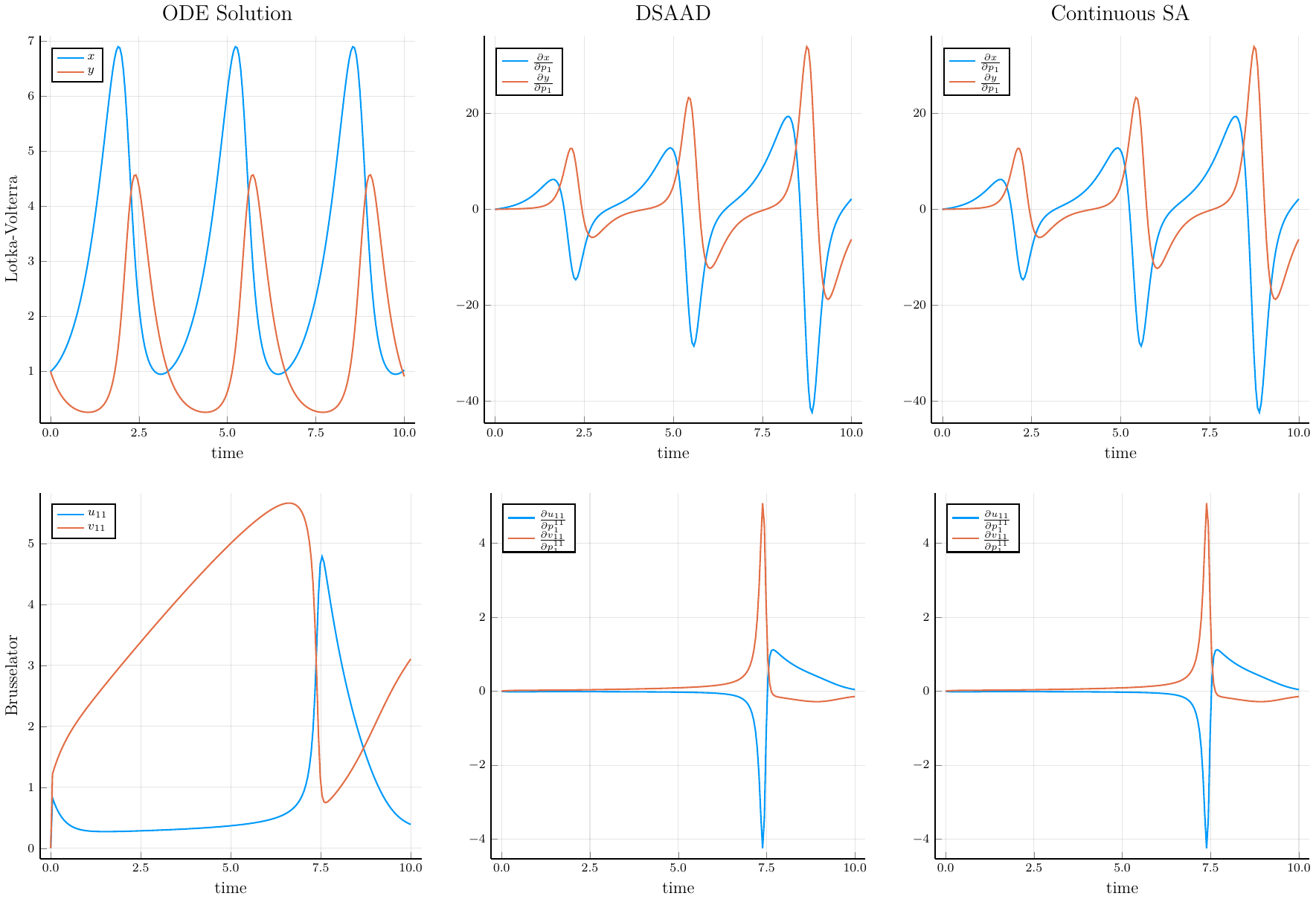}
    \end{center}
    \caption{{\bf Model Sensitivities for DSAAD and CSA.} \textbf{Top row:} For the Lotka-Volterra model, the results are shown for $t \in [0,10]$. \textbf{Bottom row:} The Brusselator PDE is discretized using a $3 \times 3$ uniform grid on the domain $[0,1]\times[0,1]$, and the resulting ODE system is solved for $t \in [0,10]$. The results for $u$ and $v$ on the $(1,1)$ grid point are shown.}
    \label{fig:ode_sens}
\end{figure}

\section{Benchmark Results\label{sec:Discrete-and-Continuous}}

\subsection{Forward-Mode Sensitivity Performance Comparisons\label{subsec:Forward-Mode-Sensitivity}}

To test the relative performance of discrete and continuous sensitivity
analysis, we utilized packages from the Julia programming language.
The method for continuous
sensitivity analysis is implemented in the DiffEqSensitivity.jl package
by directly extending a user-given ordinary differential equation.
This is done by extending the initial condition vector and defining
a new derivative function $\tilde{f}$ which performs $f$ on the
first $N$ components and adds the sensitivity equations. For performance,
construction of the Jacobian can be avoided by utilizing AD for vector-Jacobian and Jacobian-vector products. By seeding the Dual numbers to have partials $v$ and applying forward applications of $f$, the resulting
output is the desired $\frac{\partial f}{\partial u}v$. Similarly, seeding
on Tracked reals for reverse-mode autodifferentiation results in
$v' \frac{\partial f}{\partial u}$ which is the other desired quantity.
As a comparison, full Jacobian implementations were also explored, either
via a user-given analytical solution or automatic differentiation.
We also include numerical differentiation performed
by the FiniteDiff.jl package.

The test problems were solved with the various sensitivity analysis methods
and the timings are given in Table \ref{tab:Forward-Sensitivity-Analysis-Bench}.
In all of these benchmarks DSAAD performs well, being the fastest or nearly the fastest
in all cases. While continuous sensitivity analysis does not necessarily require building the
Jacobian, inspection of the generated LLVM from the forward mode reveals that the compiler
is able to fuse more operations between the Jacobian-vector product and the other parts of
the calculation, effectively improving common subexpression elimination (CSE) at the compiler
level further than implementations which call functions for the separate parts of the calculation.
Thus, given the equivalence of discrete forward sensitivities and forward-mode AD, these results
showcase that sufficiently optimized forward-mode AD methods will be preferable in most cases.

\begin{center}
\begin{table}
\begin{centering}
{\tiny
\begin{tabular}{|c|c|c|c|c|}
\hline
  Method/Runtime & LV ($\mu$s) & BRUSS (s) & POLLU (ms) &
  PKPD (ms)\tabularnewline
  \hline
  \hline
  DSAAD &                       174  & 1.94 & 12.2  & 2.66 \tabularnewline
  \hline
  CSA User-Jacobian &         429 & 727  & 572 & 17.3 \tabularnewline
  \hline
  CSA AD-Jacobian &           998 & 168 & 629 & 13.7 \tabularnewline
  \hline
  CSA AD-$Jv$ seeding &       881 & 189  & 508 & 8.46 \tabularnewline
  \hline
  Numerical Differentiation & 807 & 1.58 & 24.9 & 17.4 \tabularnewline
\hline
\end{tabular}}
\par\end{centering}
\caption{\textbf{Forward Sensitivity Analysis Performance Benchmarks.} The
  Lotka-Volterra model used by the benchmarks is the same as in Figure
  \ref{fig:ode_sens}, while the Brusselator benchmarks use a finer $5 \times 5$
  grid with the same solution domain, initial values and parameters. The
  \texttt{Tsit5} integrator is used for the Lotka-Volterra and the PKPD model.
  The \texttt{Rodas5} integrator is used for the Brusselator and POLLU.}
\label{tab:Forward-Sensitivity-Analysis-Bench}
\end{table}
\par\end{center}

\subsection{Adjoint Sensitivity Performance Comparisons\label{subsec:Adjoint-Sensitivity}}

For adjoint sensitivity analysis, adjoint DSAAD
programs were produced using a combination of the generic
DifferentialEquations.jl integrator with the tape-based automatic
differentiation implementation of ReverseDiff.jl. DiffEqSensitivity.jl
provides an implementation of CASA which
saves a continuous solution for the forward pass of the solution and utilizes
its interpolant in order to calculate the requisite Jacobian and gradients
for the backwards pass. While this method is less memory efficient than
checkpointing or re-solving schemes \cite{alan_c._hindmarsh_sundials:_2005},
it only requires a single forward numerical solution and thus was demonstrated
as more runtime optimized for sufficiently small models. Thus while
DifferentialEquations.jl contains a checkpointed adjoint implementation, checkpointing was turned off for the interpolating and backsolve schemes to allow them as much performance as possible.

The timing results of these methods on the test problems are given in
Table \ref{tab:Adjoint-Sensitivity-Analysis-Bench}. 
These results show a clear performance advantage for forward-mode DSAAD
over the other choices on sufficiently small models.

\begin{center}
\begin{table}
\begin{centering}
{\tiny
\begin{tabular}{|c|c|c|c|c|}
\hline
  Method/Runtime & LV  ($\mu$s) & BRUSS (s) & POLLU (s) &
  PKPD (ms)\tabularnewline
  \hline
  \hline
  Forward-Mode DSAAD &         279  & 1.80 & 0.010  & 	5.81 \tabularnewline
  \hline
  Reverse-Mode DSAAD &         5670 & 19.1  & 0.194 & 133 \tabularnewline
  \hline
  CASA User-Jacobian (interpolating) &         549 & 25.1  & 9.18 & 	6.48 \tabularnewline
  \hline
  CASA AD-Jacobian (interpolating) &           636 & 11.8 & 16.1 & 	5.23 \tabularnewline
  \hline
  CASA AD-$v'J$ seeding (interpolating) &      517 & 1.59  & 2.12 & 2.13 \tabularnewline
  \hline
  \hline
  CASA User-Jacobian (quadrature) &     693 & 0.964  & 1.82 & 	4.88 \tabularnewline
  \hline
  CASA AD-Jacobian (quadrature) &           825 & 2.17 & 6.19 & 	4.97 \tabularnewline
  \hline
  CASA AD-$v'J$ seeding (quadrature) &      707 & 0.461  & 1.30 & 2.94 \tabularnewline
  \hline
  CASA User-Jacobian (backsolve) &     813 & N/A  & N/A & 	N/A \tabularnewline
  \hline
  CASA AD-Jacobian (backsolve) &           941 & N/A & N/A & 	N/A \tabularnewline
  \hline
  CASA AD-$v'J$ seeding (backsolve) &      760 & N/A  & N/A & N/A \tabularnewline
  \hline
  Numerical Differentiation & 811 & 2.48 & 0.044  & 	20.8 \tabularnewline
\hline
\end{tabular}}
\par\end{centering}
\caption{\textbf{Adjoint Sensitivity Analysis Performance Benchmarks.} The Lotka-Volterra and Brusselator models used by the benchmarks are the same as in Figure \ref{fig:ode_sens}.  The integrators used for the benchmarks are: \texttt{Rodas5} for Brusselator and POLLU, and \texttt{Tsit5} for Lotka-Volterra and PKPD.}
\label{tab:Adjoint-Sensitivity-Analysis-Bench}
\end{table}
\par\end{center}

\subsection{Adjoint Sensitivity Scaling}
\begin{figure}
    \begin{center}
        \includegraphics[width=\linewidth]{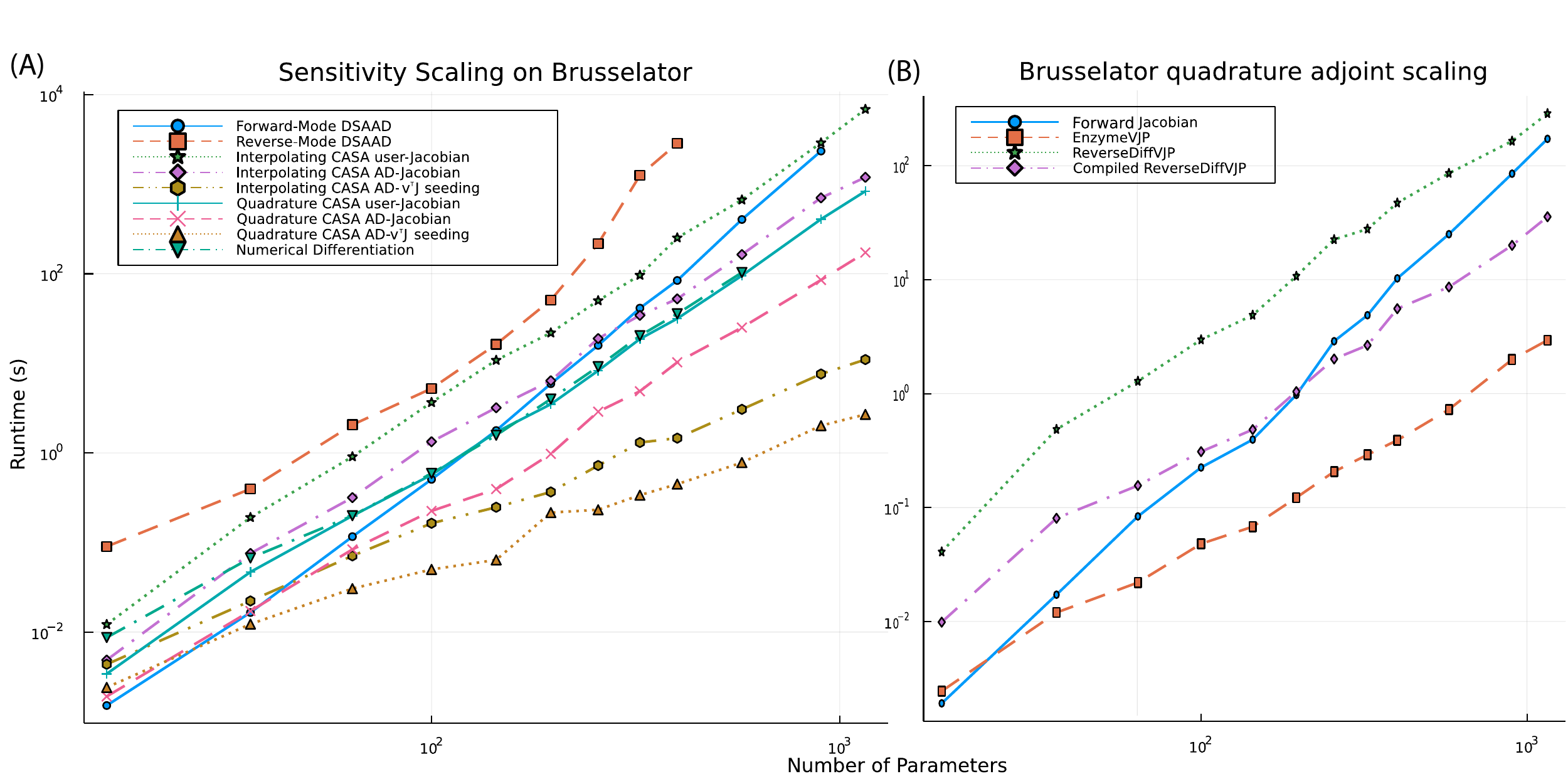}
    \end{center}
    \caption{\textbf{Brusselator Scaling Benchmarks.} A: The Brusselator problem 
    was solved with varying dimension $N$ with each of the different 
    sensitivity analysis methods to determine their scalability with respect
    to number of parameters. Depicted on the x-axis is the number of parameters
    $4N^2$ in log-scale with the y-axis being the runtime in seconds in log-scale. All
    $v'J$ seeding uses Enzyme. B: The same Brusselator, comparing between reverse-mode
    AD techiques used for $v'J$ seeding.}
    \label{fig:bruss_scaling}
\end{figure}
The previous tests all showed that on small models forward-mode via AD
was advantageous to reverse-mode and adjoint methods. However, the advantage
of adjoint methods comes in their ability to scale, with additive scaling with respect to the number of ODEs and parameters instead of multiplicative.
Thus we decided to test the scaling of the methods on the Brusselator partial differential equation.
For an $N\times N$ discretization in space, this problem has $2N^2$ ODE terms 
and  $4N^2$ parameters. The timing results for the adjoint methods and
forward-mode DSAAD are shown in Figure \ref{fig:bruss_scaling}. 
Figure \ref{fig:bruss_scaling}A demonstrates that as
$N$ increases, there is a point at which CASA becomes more efficient than
DSAAD. This cutoff point seems to be around 50 ODEs+parameters when $v'J$ seeding is used and around 100 ODEs+parameters when $v'J$ seeding is not used with quadrature adjoints, to around 150 ODEs+parameters without $v'J$ seeding with interpolating adjoints (noted in the discussion section as the method equivalent to SUNDIALS \cite{alan_c._hindmarsh_sundials:_2005}). This identifies a ballpark range of around 100 combined ODEs+parameters at which practitioners should consider changing from forward sensitivity approaches to adjoint methods.

The most efficient CASA method is the quadrature-based method. This aligns with the theoretical understanding of the method. Stiff ODE solvers scale like $\mathcal{O}(n^3)$ where $n$ is the number of ODEs. For the adjoint calculation, the quadrature based method has $n=2N^2$, the number of ODEs of the forward pass. For the interpolating adjoint method, $n=6N^2$, or the number of parameters plus the number of ODEs, which is then amplified by the cubic scaling of the linear solving. We note in passing that the backsolve technique requires $n=8N^2$, or double the number of ODEs plus the number of parameters, but is unstable on stiff equations and thus is not effective on this type of problem (indeed, on these benchmarks it was unstable and thus omitted). 
We note that the inability for the reverse-mode DSAAD to scale comes
from the tape generation performance on nonlinear models. The Tracker types in the differential
equation solvers have to utilize scalar tracking instead of tracking array
primitives, greatly increasing the size of the computational graph and the
memory burden. When array primitives are used, mutation is not allowed which
decreases the efficiency of the solver more than the gained efficiency. This is a
general phenomena of tape-based automatic differentiation methods and has been
similarly noted in PyTorch \cite{NEURIPS2019_9015} and TensorFlow Eager \cite{agrawal2019tensorflow}, demonstrating that the effective
handling of this approach would require alternative AD architectures.
These issues may be addressed in the next generation reverse-mode source-to-source 
AD packages like Zygote \cite{innes_dont_2018} or Enzyme \cite{enzymeNeurips} by not relying on tape 
generation.

Figure \ref{fig:bruss_scaling}B further highlights the importance of the reverse-mode
architecture by comparing between different reverse-mode choices used for the $v'J$
calculation internal to the CASA. By default, ReverseDiff.jl is a tape-based AD 
which has trouble scaling to larger problems. Compiled ReverseDiff.jl does a single
trace through the program to generate a static description of the code which it
optimizes and compiles (this is thus not compatible with code that is non-static, such
as having branches which change due to $u$ or $t$). Enzyme.jl compiles a static SSA-form
LLVM representation of the user's $f$ function and uses compiler passes to generate the
$v'J$ calculation. By acting at the lowest level, it can generate code
after other code optimizations have been applied, which leads to noticeable performance gains
of around an order of magnitude on this nonlinear partial differential equation application.

\subsection{Parameter Estimation Performance Comparisons\label{subsec:Application:-Parameter-Estimatio}}

To test the effect of sensitivity analysis timings in applications, we
benchmarked parameter estimation performed with an $L^2$ loss function on
generated data for each of the models. The data was generated by solving the ODE with the parameters defined in the Appendix and sampling at evenly spaced time points (100 points for Lotka-Volterra, 20 for Brusselator, 10 for POLLU, and 41 for PK/PD). Each of the parameter estimations were done using the  BFGS local optimizer from Optim.jl \cite{mogensen_optim:_2018}
and were ran until the optimizer converged to the optima to a tolerance of $10^{-6}$ (results were checked
for proper convergence). Each method started
from the same initial condition which was a perturbation of the true parameters. For Lotka-Volerra, the initial parameter values were $\frac{4}{5}$ the true values, for PK/PD $0.95p_0+0.001$, and for the other models $\frac{9}{10}$. The timings are shown in Table \ref{tab:Parameter-Estimation-Bench}.
While not as pronounced as the pure sensitivity calculations, these benchmarks
show that utilizing a more efficient sensitivity analysis calculation does
give a performance advantage in the application. Additional performance
disadvantages for numerical differentiation could be attributed to the
increased numerical error in the gradient which notably caused more iterations
for the optimizer.

\begin{center}
\begin{table}
\begin{centering}
{\tiny
\begin{tabular}{|c|c|c|c|c|}
\hline
  Method/Runtime & LV ($s$) & BRUSS ($s$) & POLLU ($s$) &
  PKPD ($s$)\tabularnewline
  \hline
  \hline
  Forward-Mode DSAAD &        3.51  & 0.760 & 0.025  & 9.54 \tabularnewline
  \hline
  CSA User-Jacobian &         1.47 &  52.2  & 0.931 & 	12.3 \tabularnewline
  \hline
  CSA AD-Jacobian &           1.82 &  51.9 & 0.932 &  12.5 \tabularnewline
  \hline
  CSA AD-$Jv$ seeding &       1.84 &  51.9  & 0.928 & 12.5 \tabularnewline
  \hline
  Reverse-Mode DSAAD &        5.87 &  14.9  & 0.560 & 238 \tabularnewline
  \hline
  CASA User-Jacobian (interpolating) &     0.030 & 12.4 & 0.238 & 42.5 \tabularnewline
  \hline
  CASA AD-Jacobian (interpolating) &       0.031 & 5.58 & 0.414 & 25.0 \tabularnewline
  \hline
  CASA AD-$v'J$ seeding (interpolating) &  0.028 & 2.26 & 0.060 & 13.4 \tabularnewline
  \hline
  CASA User-Jacobian (quadrature) &        0.122 & 1.37 & 0.053 & 21.2 \tabularnewline
  \hline
  CASA AD-Jacobian (quadrature) &          0.078 & 1.67 & 0.144 & 6.35 \tabularnewline
  \hline
  CASA AD-$v'J$ seeding (quadrature) &     0.065 & 1.20 & 0.038 & 18.3 \tabularnewline
  \hline
  CASA User-Jacobian (backsolve) &        2.95 & N/A & N/A & N/A \tabularnewline
  \hline
  CASA AD-Jacobian (backsolve) &          1.35 & N/A & N/A & N/A \tabularnewline
  \hline
  CASA AD-$v'J$ seeding (backsolve) &     1.30 & N/A & N/A & N/A \tabularnewline
  \hline
  Numerical Differentiation & 0.105 & 8.17 & 0.110  & 168 \tabularnewline
\hline
\end{tabular}}
\par\end{centering}
\caption{\textbf{Parameter Estimation Benchmarks.} The
  Lotka-Volterra model used by the benchmarks is the same as in Figure
  \ref{fig:ode_sens}, while the Brusselator benchmarks use a finer $5 \times 5$
  grid with the same solution domain, initial values and parameters. The \texttt{Tsit5} integrator is used for Lotka-Volterra and PKPD, and the \texttt{Rodas5} integrator is used for Brusselator and pollution.}
\label{tab:Parameter-Estimation-Bench}
\end{table}
\par\end{center}

\section{DSAAD Generalizes To Hybrid, Delay, and Differential-Algebraic Differential Equations\label{sec:Generalizes}}

We compared the flexibility of sensitivity analysis approaches in
order to understand their relative merits for use in a general-purpose
differential equation package. First we analyze the ability for the
sensitivity analysis approaches to work with event handling. Event-handling
is a feature of differential equation packages which allows users
to provide a rootfinding function $g(u,p,t)$ at which a discontinuity
(of the user's choice) is applied at every time point where $g(u,p,t)=0$
(such equations are also known as hybrid differential equations).
Automatic differentiation approaches to discrete sensitivity analysis
directly generalize to handling this case by propagating Dual numbers
through the event handling code. On the other hand, continuous sensitivity
analysis approaches can require special handling in order to achieve
correctness. There are two ways which sensitivities will not propagate:
\begin{enumerate}
\item Standard continuous sensitivity analysis does not take into account
the sensitivity of the time point of the discontinuity to the parameters.
\item Standard continuous sensitivity analysis does not take into account
the possibility of the discontinuity's amount being parameter-dependent.
\end{enumerate}
These points can be illustrated using a single state linear control problem where $x$ is the signal responsible for the control of $y$. This results in
a first-order linear hybrid ordinary differential equation system:
\begin{align}
    \frac{dx}{dt} & = -a,\nonumber\\
    \frac{dy}{dt} & = b,\label{eq:event}
\end{align}
where $a, b > 0$ are parameters, and the rootfinding function is $g(x,y,p,t)=x$.
At zero-crossings the parameter $b$ is set to $0$, effectively turning
off the second equation.
For the initial condition $(x(0), y(0)) = (1,0)$, the system's analytical
solution is
\begin{align}
    x(t) &= 1 - at,\nonumber\\
    y(t) &= \begin{cases}
        bt & (t < t^*)\\
        bt^* = \frac{b}{a} & (t \ge t^*),
    \end{cases} \label{eq:event-analytical}
\end{align}
where $t^* = 1/a$ is the crossing time, which depends on the parameters.
Furthermore, the amount of jump for the discontinuity of $dy/dt$ is also
parameter dependent.

\begin{table}
\begin{centering}
{\tiny
\begin{tabular}{|c|c|c|c|c|}
\hline
Method & $\partial x(1)/\partial a$ & $\partial y(1)/\partial a$ & $\partial x(1)/\partial b$ & $\partial y(1)/\partial b$\tabularnewline
\hline
\hline
Analytical Solution & -1.0 & -0.25 & 0 & 0.5 \tabularnewline
\hline
DSAAD & -1.0 & -0.25 & 5.50e-11 & 0.5 \tabularnewline
\hline
CSA & -1.0 & 0.0 & 0.0 & 1.0 \tabularnewline
\hline
\end{tabular}}
\par\end{centering}
\caption{{\bf Sensitivity Analysis with Events}. Shown are the results of the control problem given by Equation \ref{eq:event} with $a = 2$ and $b = 1$. It was solved on $t\in[0,1]$ with initial condition $(x(0), y(0)) = (1,0)$. The sensitivities of the two state variables are given at time $t=1$ with respect to the two parameters. The analytical solution is derived by taking derivatives directly on Equation \ref{eq:event-analytical}.}
\label{tab:event-results}
\end{table}

The sensitivity analysis results are compared to the true derivative
at time $t=1$ utilizing the analytical solution of the system
(Equation \ref{eq:event-analytical}) in Table
\ref{tab:event-results}. These results show that continuous sensitivity
analysis as defined in Equation \ref{eq:clsa} does not properly propagate
the sensitivities due to discontinuities and this results in incorrect
derivative calculations for hybrid ODE systems. Extensions to continuous
sensitivity analysis are required in order to correct these errors \cite{kirches_numerical_2006}. These have been implemented in DiffEqSensitivity.jl, though we remark from experience that such corrections are very non-trivial to implement and thus demonstrate an ease of implementation advantage for DSAAD.

Additionally, the continuous sensitivity analysis equations defined
in Equation \ref{eq:clsa} only apply to ordinary differential equations.
It has been shown that a different set of equations is required for
delay differential equations (DDEs) \cite{bredies_generalized_2013}
\begin{equation}
  \frac{d}{dt}\frac{\partial u(t)}{\partial p_i} = \frac{\partial G}{\partial u}
  \frac{\partial u}{\partial p_i}(t) + \frac{\partial G}{\partial
  \tilde{u}}\frac{\partial u}{\partial p_i}(t-\tau) +
  \frac{\partial G}{\partial p_i}(t),
\end{equation}
where $\frac{du(t)}{dt}=G(u,\tilde{u},p,t)$ is the DDE system with a single fixed time delay $\tau$,
and differential-algebraic equations (DAEs) \cite{alan_c._hindmarsh_sundials:_2005}
\begin{equation}
  \frac{\partial F}{\partial u}\frac{\partial u}{\partial p_i} + \frac{\partial
  F}{\partial \dot{u}}\frac{\partial \dot{u}}{\partial p_i} + \frac{\partial F}{\partial
  p_i} = 0,
\end{equation}
where $F(\dot{u},u,p,t)=0$ is the DAE system.
On the other hand, the discrete sensitivity analysis approach implemented
via automatic differentiation is not specialized to ordinary differentiation
equations since it automatically generates the sensitivity propagation
at the compiler-level utilizing the atomic operations inside the numerical
integration scheme. Thus these types of equations, and expanded forms such 
as DDEs with state-dependent delays or hybrid DAEs, are automatically 
supported by the connection between DifferentialEquations.jl and 
Julia-based automatic differentiation packages. Tests on the DDE and DAE
solvers confirm this to be the case.

\section{Discussion\label{sec:Discussion}}

Performant and correct sensitivity analysis is crucial to many applications
of differential equation models. Here we analyzed both the performance
and generalizability of the full Cartesian product of approaches. Our results show a strong
performance advantage for automatic differentiation based discrete
sensitivity analysis for forward-mode sensitivity analysis on sufficiently
small systems, and an
advantage for continuous adjoint sensitivity analysis for sufficiently large systems.
Notably, pure tape-based reverse-mode automatic differentiation did not perform or
scale well in these benchmarks. The implementations have been generally optimized for the usage in machine learning models which make extensive use of
large linear algebra like matrix multiplications, which decrease the
size of the tape with respect to the amount of work performed. Since
differential equations tend to be defined by nonlinear functions with
scalar operations, the tape handling to work ratio thus decreases and
is no longer competitive with other forms of derivative calculations. This technical detail affects common frameworks such as ReverseDiff.jl, PyTorch \cite{NEURIPS2019_9015}, and TensorFlow Eager \cite{agrawal2019tensorflow}, which leads to the recommendation of CASA for most frameworks. In addition, some frameworks, such as Jax \cite{jax2018github}, cannot JIT optimize the non-static computation graphs of a full ODE solver, which further leads to performance improvements of CASA in current reverse-mode implementations.
Future work should investigate adjoint DSAAD with
approaches which allow for low overhead compilation of the reverse path, such
as that of Zygote \cite{innes_dont_2018}, for which the static handling of control flow may be what is needed for this application.

One major result to note is that, in many cases of interest, runtime
overhead of the adjoint methods can be larger than the theoretical
scaling advantages. Forward-mode automatic differentiation does
not exhibit the best scaling properties but on small systems of ODEs with small numbers of parameters, both stiff and non-stiff, the forward mode methods
benchmarks as advantageous. On problems like PDEs, the scaling
of adjoint methods does confer an advantage to them when the problem
is sufficiently large but also a disadvantage when the problem
is small. Additionally, the ability to seed automatic differentiation
for the Jacobian vector multiplications was demonstrated as
very advantageous on the larger problems, showing that integrating the implementation of sensitivity analysis with automatic differentiation tools is necessary for achieving the utmost efficiency. In addition, the 
choice of automatic differentiation for the $v'J$ calculation was shown to be a major factor in performance, with the static Enzyme-based implementation giving a two order of magnitude performance advantage over the tape-based implementation even when confined to this one calculation.
Together, these results show that having the ability to choose between these different methods is essential for a software wishing to support these distinct use cases.

Future research should look into the comparative performance of DiffEqSensitivity.jl \cite{rackauckas2020universal} with the native adjoint techniques of SUNDIALS \cite{alan_c._hindmarsh_sundials:_2005} and PETSc TS \cite{abhyankar2018petsc} to see how much this result generalizes. The results here would suggest the performance difference may be around three orders of magnitude since SUNDIALS uses a method similar to the interpolating CASA with (numerical) forward-mode Jacobians. A follow-up study on multiple PDEs comparing between the software suites could be very useful to practitioners. 

While runtime performance of these methods is usually of interest,
it is important to note the memory scaling for the various adjoint
sensitivity analysis implementations. The quadrature CASA implementation 
benchmarks as the fastest for stiff ODEs, but uses a continuous solution
to the original ODE in order to generate the adjoint Jacobian and
gradients on-demand. This setup only requires a single forward ODE
solve but makes a tradeoff due to the high memory requirement to save
the full timeseries solution and its interpolating function. For example,
with the chosen 9th order explicit Runge-Kutta method due to Verner,
the total memory cost is $26NM$ since 26 internal derivative calculations
of size $N$ are utilized to construct the interpolant where $M$ is the
number of time points. The DSAAD reverse-mode AD approaches require constructing a tape for the entire
initial forward solution and thus has memory scaling similar to the
quadrature CASA, though in theory checkpointing AD systems AD could alleviate these memory issues. 
In contrast, checkpointing-based adjoint sensitivity
analysis implementations \cite{alan_c._hindmarsh_sundials:_2005} 
re-solve the forward ODE from a saved
time point (a checkpoint) in order to get the $u$ value required
for the Jacobian and gradient calculation, increasing the runtime
cost by at most 2x while decreasing the memory cost to NC where 
C is the number of checkpoints. These results show that the DiffEqSensitivity.jl 
interpolating checkpointing CASA approach may have around a 5x performance
deficit over the more optimal quadrature CASA, which is can be a useful price
to pay if memory concerns are the largest factor. The backsolve CASA approach,
while the most memory efficient, was shown to be too unstable to be useful in 
most usecases and is thus only recommended if the dynamical system is known 
to be non-stiff in advance.

These results show many advantages for the AD-based discrete sensitivity
analysis for small systems. However, there are significant engineering challenges to
the development of such integration schemes. In order for this methodology
to exist, a general function automatic differentiation tool must exist
for the programming language and the entire ODE solver must be compatible
with the AD software. This works for the Julia-based DifferentialEquations.jl
software since it contains a large number of high performance native-Julia
solver implementations. However, many other solver ecosystems do not
have such possibilities. For example, common open source packages
for solving ordinary differential equation systems include deSolve
in R \cite{soetaert_solving_2010} and SciPy for Python \cite{jones_scipy:_2001}. While general automatic differentiation
tools exist for these languages (\cite{pav_madness:_nodate}, \cite{maclaurin_autograd:_nodate}), both of these package
call out to Fortran-based ODE solvers such as LSODA \cite{petzold_automatic_nodate}, and thus
AD cannot be directly applied to the solver calls. This means DSAAD techniques may stay limited to specific software due to technical engineering issues.

When considering the maintenance of large software ecosystems,
the DSAAD approach gives many advantages if planned from the start.
For one, almost no additional code was required to be written by the
differential equation community in order for this implementation to
exist since it directly works via code generation at compile-time
on the generic functions of DifferentialEquations.jl. But an additional
advantage is that this same technique applies to the native hybrid,
delay, and differential-algebraic integrators present in the library.
DifferentialEquations.jl also allows for many other actions to occur
in the events. For example, the user can change the number of ODEs
during an event, and events can change solver internals like the current
integration time. Continuous sensitivity analysis require a
separate implementation for each of these cases with could be
costly to developer time. For DiffEqSensitivity.jl, the performant and correct
implementation of corrections for hybrid differential equations and differential-algebraic
equations took around 6 months of split between two developers (part time). This
ease of implementation aspect may be worth considering for smaller organizations.

\section*{Acknowledgments}

This work is supported by Center for Translational Medicine,
University of Maryland Baltimore School of Pharmacy. CR is partially supported
by the NSF grant DMS1763272 and a grant from the Simons Foundation
(594598, QN). Additionally, we thank the Google Summer of Code program for helping support the open source packages showcased in this manuscript.
The information, data, or work presented herein was funded in part by the Advanced Research Projects Agency-Energy (ARPA-E), U.S. Department of Energy, under Award Number DE-AR0001222 and NSF grants OAC-1835443 and IIP-1938400.



\bibliographystyle{IEEEtran}
\bibliography{IEEEabrv,bib}

\section{Appendix}
\scriptsize

\subsection{Models}

The first test problem is LV, the non-stiff Lotka-Volterra model
{\scriptsize
\begin{align}
  \frac{dx}{dt} &= p_{1}x - p_{2}xy,\nonumber \\
  \frac{dy}{dt} &= -p_{3}y + xy.\label{eq:Lotka-Volterra}
\end{align}}

with initial condition $[1.0,1.0]$ and $p = [1.5,1.0,3.0]$ \cite{murray_mathematical_2002}. 
The second model, BRUSS, is the two dimensional ($N\times N$) Brusselator stiff reaction-diffusion PDE:

{\scriptsize
\begin{align}
  \frac{\partial u}{\partial t} &= p_{2} + u^2v - (p_{1}+1)u + p_{3}(\frac{\partial^2 u}{\partial x^2}+\frac{\partial^2 u}{\partial y^2}) + f(x, y, t),\nonumber \\
  \frac{\partial v}{\partial t} &= p_{1}u - u^2v + p_{4}(\frac{\partial^2 u}{\partial x^2}+\frac{\partial^2 u}{\partial y^2}),\label{eq:Brusselator}
\end{align}}
where
{\scriptsize
\begin{align}
  f(x, y, t) = \begin{cases}
    5 & \text{ if } (x-0.3)^2+(y-0.6)^2\le 0.1^2 \text{ and } t\ge 1.1 \\
    0 & \text{ else},
  \end{cases}
\end{align}
}
with no-flux boundary conditions and $u(0,x,y) = 22(y(1-y))^{3/2}$ with
$v(0,x,y) = 27(x(1-x))^{3/2}$ \cite{ernst_hairer_solving_1991}. 
This PDE is discretized to a 
set of $N \times N \times 2$ ODEs using the finite difference
method. The parameters are spatially-dependent, $p_i=p_i(x,y)$,
making each discretized $p_i$ a $N \times N$ set of values at each
discretization point, giving a total of $4N^2$ parameters. The
initial parameter values were the uniform $p_i(x,y) = [3.4,1.0,10.0,10.0]$

% https://archimede.dm.uniba.it/~testset/report/pollu.pdf
The third model, POLLU, simulates air pollution. It is a stiff non-linear ODE system
which consists $20$ ODEs:
{\allowdisplaybreaks
{\scriptsize
\begin{align}
  \frac{du_1}{dt}  &= -p_{1}
  u_{1}-p_{10}u_{11}u_{1}-p_{14}u_{1}u_{6}-p_{23}u_{1}u_{4}- \nonumber\\
  & p_{24}u_{19}u_{1}+p_{2} u_{2}u_{4}+p_{3} u_{5}u_{2}+p_{9}u_{11}u_{2}+ \nonumber\\
  & p_{11}u_{13}+p_{12}u_{10}u_{2}+p_{22}u_{19}+p_{25}u_{20} \nonumber\\
  \frac{du_2}{dt}  &= -p_{2} u_{2}u_{4}-p_{3} u_{5}u_{2}-p_{9}
  u_{11}u_{2}-p_{12}u_{10}u_{2}+p_{1} u_{1}+p_{21}u_{19} \nonumber\\
  \frac{du_3}{dt}  &= -p_{15}u_{3}+p_{1}
  u_{1}+p_{17}u_{4}+p_{19}u_{16}+p_{22}u_{19} \nonumber\\
  \frac{du_4}{dt}  &= -p_{2}
  u_{2}u_{4}-p_{16}u_{4}-p_{17}u_{4}-p_{23}u_{1}u_{4}+p_{15}u_{3} \nonumber\\
  \frac{du_5}{dt}  &= -p_{3} u_{5}u_{2}+p_{4} u_{7}+p_{4} u_{7}+p_{6}
  u_{7}u_{6}+p_{7} u_{9}+p_{13}u_{14}+p_{20}u_{17}u_{6} \nonumber\\
  \frac{du_6}{dt}  &= -p_{6} u_{7}u_{6}-p_{8}
  u_{9}u_{6}-p_{14}u_{1}u_{6}-p_{20}u_{17}u_{6}+p_{3}
  u_{5}u_{2}+p_{18}u_{16}+p_{18}u_{16} \nonumber\\
  \frac{du_7}{dt}  &= -p_{4} u_{7}-p_{5} u_{7}-p_{6} u_{7}u_{6}+p_{13}u_{14} \nonumber\\
  \frac{du_8}{dt}  &= p_{4} u_{7}+p_{5} u_{7}+p_{6} u_{7}u_{6}+p_{7} u_{9} \nonumber\\
  \frac{du_9}{dt}  &= -p_{7} u_{9}-p_{8} u_{9}u_{6} \nonumber\\
  \frac{du_{10}}{dt} &= -p_{12}u_{10}u_{2}+p_{7} u_{9}+p_{9} u_{11}u_{2} \nonumber\\
  \frac{du_{11}}{dt} &= -p_{9} u_{11}u_{2}-p_{10}u_{11}u_{1}+p_{8} u_{9}u_{6}+p_{11}u_{13} \nonumber\\
  \frac{du_{12}}{dt} &= p_{9} u_{11}u_{2} \nonumber\\
  \frac{du_{13}}{dt} &= -p_{11}u_{13}+p_{10}u_{11}u_{1} \nonumber\\
  \frac{du_{14}}{dt} &= -p_{13}u_{14}+p_{12}u_{10}u_{2} \nonumber\\
  \frac{du_{15}}{dt} &= p_{14}u_{1}u_{6} \nonumber\\
  \frac{du_{16}}{dt} &= -p_{18}u_{16}-p_{19}u_{16}+p_{16}u_{4} \nonumber\\
  \frac{du_{17}}{dt} &= -p_{20}u_{17}u_{6} \nonumber\\
  \frac{du_{18}}{dt} &= p_{20}u_{17}u_{6} \nonumber\\
  \frac{du_{19}}{dt} &= -p_{21}u_{19}-p_{22}u_{19}-p_{24}u_{19}u_{1}+p_{23}u_{1}u_{4}+p_{25}u_{20} \nonumber\\
  \frac{du_{20}}{dt} &= -p_{25}u_{20}+p_{24}u_{19}u_{1}
\end{align}
}}
with the initial condition of $u_0 = [0, 0.2, 0, 0.04, 0, 0, 0.1, 0.3, 0.01, 0,
0, 0, 0 ,0, 0, 0, 0.007, 0, 0, 0]^T$ and parameters
$[.35, 26.6, 12,300, .00086, .00082, 15,000, .00013, 24,000, 16,500, 9,000,\\ 
.022, 12,000, 1.88, 16,300, 4,800,000, .00035,.0175, 10^9, .444e12, 1,240,\\ 
2.1, 5.78, .0474, 1,780, 3.12]$ \cite{ernst_hairer_solving_1991}.

The fourth model is a non-stiff pharmacokinetic/pharmacodynamic model (PKPD) \cite{upton_basic_2014}, which is in
the form of
{\scriptsize
\begin{align}
  \frac{dDepot}{dt}  &= -k_a Depot \nonumber\\
  \frac{dCent}{dt}   &= k_a Depot+ \nonumber\\
  & (CL+V_{max}/(K_m+(Cent/V_c))+Q_1)(Cent/V_c) + \nonumber\\
  & Q_1(Periph_1/V_{p_1}) - \nonumber\\
  & Q_2(Cent/V_c)  + Q_2(Periph_2/V_{p_2}) \nonumber\\
  \frac{dPeriph_1}{dt}   &= Q_1(Cent/V_c)  - Q_1(Periph_1/V_{p_1}) \nonumber\\
  \frac{dPeriph_2}{dt}   &= Q_2(Cent/V_c)  - Q_2(Periph_2/V_{p_2}) \nonumber\\
  \frac{dResp}{dt}   &=
  k_{in}(1-(I_{max}(Cent/V_c)^\gamma/(IC_{50}^\gamma+(Cent/V_c)^\gamma)))
  - k_{out} Resp.
  %Ev'    &= -Ka_1\cdotEv \nonumber\\
  %Cent'   &=  Ka_1\cdotEv -
  %(CL+V_{max}/(K_m+(Cent/V_c))+Q)\cdot(Cent/V_c)  +
  %Q\cdot(Periph/Vp) - \nonumber\\
  %&Q_2\cdot(Cent/V_c)  + Q_2\cdot(Periph_2/Vp_2) \nonumber\\
  %Periph' &=  Q\cdot(Cent/V_c)  - Q\cdot(Periph/V_p) \nonumber\\
  %Periph_2' &=  Q_2\cdot(Cent/V_c)  - Q_2\cdot(Periph_2/Vp2) \nonumber\\
  %Resp'   &=
  %Kin\cdot(1-(I_{max}\cdot(Cent/V_c)^\gamma/(IC50^\gamma+(Cent/V_c)^\gamma)))
  %- K_{out}\cdotResp.
\end{align}}
with the initial condition of $[100.0,0.0,0.0,0.0,5.0]$. $k_a = 1$ is the absorption rate of drug into the central compartment from the dosing compartment, $CL = 1$ is the clearance parameter of drug elimination, $V_c = 20$ is the central volume of distribution, $Q_1 = 2$ is the inter-compartmental clearance between central and first peripheral compartment, $Q_2 = 0.5$ is the inter-compartmental clearance between central and second peripheral compartment, $V_{p_1} = 10$ is the first peripheral compartment distribution volume, $V_{p_2} = 100$ is the second peripheral compartment distribution volume, $V_{max} = 0$ is the maximal rate of saturable elimination of drug, $K_m = 2$ is the Michaelis-Mentens constant, $k_{in} = 10$ is the input rate to the response (PD) compartment with a maximal inhibitory effect of $I_{max} = 1$, $IC_{50} = 2$ is a parameter for the concentration at 50\% of the effect and $k_{out} = 2$ is the elimination rate of the response,
and $\gamma = 1$ is the model sigmoidicity. Additional doses of $100.0$ are applied to the $Depot$ variable at every $24$ time units.

\end{document}